\documentclass[12pt]{amsart}

\usepackage{enumitem}
\usepackage{psfrag}
\usepackage{graphicx}
\usepackage{pinlabel}
\usepackage{graphicx}
\usepackage{color}
\usepackage{amsmath}
\usepackage{amssymb}
\usepackage{amscd}
\usepackage{picinpar}
\usepackage{times}
\usepackage{pb-diagram}
\usepackage{graphicx}
\usepackage{wrapfig}
\usepackage{xspace}
\usepackage{hyperref}
\usepackage{url}
\usepackage{caption}
\usepackage{subcaption}
\usepackage{fancyhdr}
\newtheorem{theorem}{Theorem}[section]
\newtheorem{lemma}[theorem]{Lemma}
\newtheorem{proposition}[theorem]{Proposition}
\newtheorem{corollary}[theorem]{Corollary}
\newtheorem{definition}[theorem]{Definition}

\newtheorem{remark}[theorem]{Remark}
\newtheorem*{theorem*}{Theorem}
\def\C{\mathbb{C}}
\def\R{\mathbb{R}}

\def\T{\mathcal{T}}

\def\D{\mathbb{D}}

\usepackage{overpic}
\newcommand{\bal}{\begin{aligned}}      \newcommand{\eal}{\end{aligned}}
\newcommand{\ba}{\begin{array}}      \newcommand{\ea}{\end{array}}
\newcommand{\bc}{\begin{center}}     \newcommand{\ec}{\end{center}}
\newcommand{\be}{\begin{enumerate}}  \newcommand{\ee}{\end{enumerate}}
\newcommand{\beq}{\begin{eqnarray}}  \newcommand{\eeq}{\end{eqnarray}}
\newcommand{\beQ}{\begin{eqnarray*}} \newcommand{\eeQ}{\end{eqnarray*}}
\newcommand{\bi}{\begin{itemize}}    \newcommand{\ei}{\end{itemize}}
\newcommand{\bt}{\begin{tabular}}    \newcommand{\et}{\end{tabular}}
\newcommand{\bdm}{\begin{displaymath}} \newcommand{\edm}{\end{displaymath}}

\begin{document}

\title{The convergence of inversive distance circle packings to the Riemann mapping}
\author{Yuxiang Chen, Yanwen Luo, Xu Xu}

\address{School of Mathematics and Statistics, Wuhan University, Wuhan, 430072, P.R.China}
 \email{chenyuxiang@whu.edu.cn}

\address{Department of Mathematics, Oklahoma State University, Stillwater, 74074, USA
 }
 \email{yanwen.luo@okstate.edu}

\address{School of Mathematics and Statistics, Wuhan University, Wuhan, 430072, P.R.China}
 \email{xuxu2@whu.edu.cn}




\begin{abstract}
In \cite{BS}, Bowers and Stephenson  introduced the notion of inversive distance circle packings as a natural generalization of Thurston's
circle packings \cite{Th}.
They further conjectured that discrete conformal maps induced by inversive distance circle packings converge to the Riemann mapping.
In this paper, we prove Bowers-Stephenson's conjecture for Jordan domains
by establishing a solvability theorem of certain prescribing combinatorial curvature problems for inversive distance circle packings.
\end{abstract}

\maketitle


\section{Introduction}
In \cite{Th2}, Thurston proposed a constructive approach to the Riemann mapping theorem by approximating conformal mappings in simply connected domains using circle packings. Thurston conjectured that the discrete conformal maps induced by circle packings converge to the Riemann mapping. Thurston's conjecture has been proved elegantly by Rodin-Sullivan \cite{RS}.
Since then, there have been lots of important works on the convergence of discrete conformal maps to the Riemann mapping. See \cite{Bucking, GLW, HS1, HS2, LSW,WZ} and others.
Motivated by Thurston's circle packings \cite{Th}, Bowers-Stephenson \cite{BS} introduced the notion of
inversive distance circle packings and conjectured that the Riemann mapping
could be approximated by inversive distance circle packings.
In this paper, we prove Bowers-Stephenson's conjecture for Jordan domains as a counterpart of Thurston's conjecture \cite{Th2} in the setting of
circle packings.
The main idea comes from the recent work of Luo-Sun-Wu \cite{LSW}.

\subsection{Definitions}
Suppose $S$ is a topological surface possibly with boundary and $\mathcal{T}$
is a triangulation of $S$.
We use $V = V(\mathcal{T})$, $E = E(\mathcal{T})$ and $F = F(\mathcal{T})$ to
denote the set of vertices, edges, and faces of $\mathcal{T}$ respectively.
A piecewise linear metric $d$ (PL metric for simplicity) on $(S, \mathcal{T})$
is a flat cone metric on $S$ such that each face in $F$ in the metric $d$ is isometric to a non-degenerate Euclidean triangle.
In this case, one can represent the PL metric on $(S, \mathcal{T})$ as a length function $l: E\rightarrow \mathbb{R}_{>0}$,
which satisfies the triangle inequality for any face in $F$.
Conversely, given a function $l: E\rightarrow \mathbb{R}_{>0}$ satisfying
the triangle inequality, one can construct a PL metric on $(S, \mathcal{T})$
by isometrically gluing Euclidean triangles along edges in pairs. Hence, we also refer to a PL metric on $(S, \mathcal{T})$
as a function $l: E\rightarrow \mathbb{R}_{>0}$ satisfying the strict triangle inequality for any face in $F$.

For a PL metric $l: E\rightarrow \mathbb{R}_{>0}$ on $(S, \mathcal{T})$, the combinatorial curvature is a map $K: V\rightarrow (-\infty, 2\pi)$
sending an interior vertex $v\in V$ to $2\pi$ minus the sum of angles at $v$ and
a boundary vertex $v\in V$ to $\pi$ minus the sum of angles at $v$.
The combinatorial curvature $K$ for a PL metric on $(S, \mathcal{T})$ satisfies the discrete Gauss-Bonnet formula
\begin{equation}\label{discrete Gauss-Bonnet formula}
  \sum_{v\in V}K(v)=2\pi\chi(S),
\end{equation}
where $\chi(S)$ is the Euler characteristic of the surface. A vertex $v$ is flat in a PL metric if $K(v) = 0$. A PL metric is flat if all interior vertices are flat.

\begin{definition}[\cite{BS}]\label{discrere conformal for idcp}
Suppose $(S, \mathcal{T})$ is a triangulated surface with a weight $I: E\rightarrow (-1, +\infty)$.
  A PL metric $l: E\rightarrow \mathbb{R}_{>0}$ on $(S, \mathcal{T})$
  is an inversive distance circle packing metric on the weighted triangulated surface $(S, \mathcal{T}, I)$
  if there exists a function $u: V \to \mathbb{R}$  such that for any edge $e\in E$ with vertices $v$ and $v'$, the length $ l(e)$ is given by
\begin{equation}
\label{length1 introduction}
  l(e)=\sqrt{e^{2u(v)}+e^{2u(v')}+2 I(e)e^{u(v)+u(v')}}.
\end{equation}
The function $u: V \to \mathbb{R}$ is called a label on $(S, \mathcal{T}, I)$.
Two inversive distance circle packing metrics $(S, \mathcal{T}, I, l)$ and $(S, \mathcal{T}, \tilde I, \tilde l)$ are conformally equivalent if $I = \tilde I$.
In this case, we set $w=\tilde{u}-u$ and denote this relation as $\tilde l =w*l$.  The function $w$ is called a discrete conformal factor on $(S, \mathcal{T}, I, l)$.
\end{definition}

If we set $r(v)=e^{u(v)}$ for $v\in V$, the weight $I(e)$ in (\ref{length1 introduction}) is the inversive distance of
the two circles centered at $v$ and $v'$ with radii $r(v)$ and $r(v')$ respectively.
The map $r: V\rightarrow (0, +\infty)$ is referred as an \textit{inversive distance circle packing} on the weighted
triangulated surface $(S, \mathcal{T}, I)$.
Thurston's circle packing \cite{Th} is a special type of
inversive distance circle packing with $I\in [0,1]$ in (\ref{length1 introduction}). An excellent source for the comprehensive theory of circle packings is
 \cite{St}.

\subsection{The main theorem}

Let $\Omega$ be a Jordan domain in the plane with three distinct boundary points $p,q,r$ specified. By the Riemann mapping theorem,
there exists a conformal map from $\Omega$ to the interior of an equilateral Euclidean triangle $\triangle ABC$ with unit edge length.
By Caratheodory's extension theorem \cite{P book},
this conformal map could be uniquely extended to be a homeomorphism $g$
from $\overline{\Omega}$ to $\triangle ABC$ with $p,q,r$ sent to $A, B,C$ respectively.
The map $g$ and $g^{-1}$ are referred as the \textit{Riemann mapping} for $(\Omega, (p,q,r))$.
Let $(D,\mathcal{T}, I)$ be an oriented weighted polygonal disk in the plane with three distinct boundary vertices $p,q,r$ and
$l$ be a flat inversive distance circle packing metric on $(D,\mathcal{T}, I)$.
Suppose that there exists a function $w: V\rightarrow \mathbb{R}$ such that $l^*=w*l$ is an inversive distance circle packing metric
on $(D,\mathcal{T}, I)$ with total area $\frac{\sqrt{3}}{4}$, combinatorial curvature $\frac{2\pi}{3}$ at $p,q,r$, and flat at other vertices.
Then $(D,\mathcal{T}, l^*)$ is isometric to
a triangulated unit equilateral triangle $(\triangle ABC, \mathcal{T}')$ with some triangulation $\mathcal{T}'$ and
the standard flat metric.
Let $f$ be the orientation-preserving piecewise linear map induced by the map sending the vertices of $\mathcal{T}$ to the corresponding vertices of $\mathcal{T}'$ such that $f(A) = p$, $f(B) = q$ and $f(C) = r$. The map
$f$ is called the \textit{discrete conformal map} associated to $(D,\mathcal{T}, I, l, \{p,q,r\})$.

We prove the following theorem on the convergence of discrete conformal maps induced by a specific sequence of inversive distance circle packings on $\Omega$, which provides an affirmative answer to Bowers-Stephenson's conjecture on the convergence inversive distance circle packings.

\begin{theorem}\label{conv introduction}
Let $\Omega$ be a Jordan domain in the complex
plane with three distinct boundary points $p,q, r$ specified. Let $f$ be the Riemann mapping from the equilateral triangle $\triangle ABC$ to $\overline{\Omega}$   such that $f(A) = p$, $f(B) = q$, $f(C) = r$.
Then there exists a sequence of weighted triangulated polygonal disks $(\Omega_{n}, \T_n,
I_n, (p_n, q_n, r_n))$ with inversive distance circle packing metrics $l_n$, where $\T_n$ is a
triangulation of $\Omega_n$, $I_n: E_n\rightarrow (1, +\infty)$ is a weight defined on $E_n = E( \mathcal{T}_n$)
and $p_n, q_n, r_n$ are three distinct boundary vertices of $\mathcal{T}_n$, such that
\begin{enumerate}
  \item[(a)] $\Omega=\cup_{n=1}^{\infty} \Omega_n$ with $\Omega_n \subset
\Omega_{n+1}$, and  $\lim_n p_n =p$, $\lim_n q_n =q$, $\lim_n r_n=r$.
  \item[(b)] discrete conformal maps $f_n$ from $\triangle ABC$ to $(\Omega_{n},
\T_n,I_n, l_n)$  with  $f_n(A) = p_n$, $f_n(B) = q_n$, $f_n(C) = r_n$  exist.
  \item[(c)] discrete conformal maps $f_n$ converge uniformly to the Riemann mapping $f$.
\end{enumerate}
\end{theorem}

In Rodin-Sullivan's convergence theorem for circle packings \cite{RS}, approximating triangulated disks can be arbitrarily selected, as the Koebe-Andreev-Thurston theorem guarantees the existence of circle packings for any triangulated disks. However, the discrete conformal map does not exist for general inversive distance circle packings with inversive distance $I:E\to (1, +\infty)$ if  the triangulation of the disk is fixed. Theorem \ref{conv introduction} resolves this issue through sufficiently fine triangulation subdivision, which produces a specific sequence of discrete conformal maps for inversive distance circle packings.   In the rest of this paper, we assume that $I:E\to (1,+\infty)$ unless otherwise stated. This condition corresponds to the ``S-packings" introduced by Bowers-Stephenson \cite{BS}.

In a recent novel work of Bobenko-Lutz \cite{BL,BL2}, by allowing the modification of triangulations, they solve the fundamental problems on the existence and uniqueness of discrete conformal metrics on closed Euclidean and hyperbolic polyhedral surfaces induced by inversive distance circle packings. Moreover, they established the beautiful relationships between the inversive distance circle packings on two dimensional polyhedral surfaces and the convex realization for convex hyperideal polyhedral cusps in three dimensional hyperbolic space. It is natural to ask whether the construction in \cite{BL} provides a family of discrete conformal maps converging to the Riemann mapping.

\subsection{Organization} The paper is organized as follows.
In Section \ref{section 2}, we give some preliminaries on inversive distance circle packings.
In Section \ref{section 3}, we solve some prescribing combinatorial curvature problem for inversive distance circle packings, establishing the existence of discrete conformal maps approximating the Riemann map. 
In Section \ref{section 4},  we prove Theorem \ref{conv introduction} based on the properties of quasiconformal maps.

\subsection{Acknowledgment}
The research of Xu Xu is supported by National Natural Science Foundation of China under grant no. 12471057.
Yanwen Luo is partially supported by Pacific Institute for the Mathematical Science, Professor Bojan Mohar from Simon Fraser University, and Professor Ryan Budney from University of Victoria. Yanwen Luo would like to thank Professor Feng Luo from Rutgers University for helpful discussions about this project.

\section{Preliminaries on inversive distance circle packings}\label{section 2}


\subsection{Admissible space of inversive distance circle packings}
Let $(S, \mathcal{T}, I)$ be a weighted triangulated surface.
 We use $v_i$ to denote a vertex in $V$, $e_{ij} = v_iv_j$ to denote an edge in $E$ and $\triangle v_iv_jv_k$ to denote a face in $F$.
We will denote $f_i = f(v_i)$ if $f$ is a function defined on $V$,  $f_{ij} = f(v_iv_j) = f(e_{ij})$ if $f$ is a function defined on $E$, and $f_{ijk} = f(\triangle v_iv_jv_k)$ if $f$ is a function defined on $F$.

For any function $u: V\rightarrow \mathbb{R}$, the formula (\ref{length1 introduction})
produces a positive function $l$ on $E$.
However, for a face $\triangle v_iv_jv_k$ in $(S, \mathcal{T}, I)$, the positive numbers $l_{ij}, l_{ik}, l_{jk}$ may not satisfy
the triangle inequality
\begin{equation}\label{strict triangle inequality}
l_{rs}< l_{rt}+l_{st}, \{r,s,t\}=\{i,j,k\}.
\end{equation}
The label $u: V\rightarrow \mathbb{R}$ is said to be \textit{admissible}
if the function $l: E\rightarrow (0, +\infty)$ determined by $u: V\rightarrow \mathbb{R}$
via the formula (\ref{length1 introduction}) satisfies
the triangle inequality (\ref{strict triangle inequality}) for every face in $(S, \mathcal{T}, I)$.
We also say that the corresponding inversive distance circle packing $r:V\rightarrow \mathbb{R}_{>0}$
on $(S, \mathcal{T}, I)$ with $r_i=e^{u_i}$ is admissible, if it causes no confusion in the context.
The admissible space of inversive distance circle packings on $(S, \mathcal{T}, I)$ consists of all the admissible inversive distance circle packings on
$(S, \mathcal{T}, I)$.
For an admissible inversive distance circle packing $r$ on $(S, \mathcal{T}, I)$, every face in $(S, \mathcal{T}, I)$ is isometric to a \textit{non-degenerate} Euclidean triangle
with edge lengths given by (\ref{length1 introduction}).
We also say that $r: V\rightarrow \mathbb{R}_{>0}$ generates a PL metric on $(S, \mathcal{T}, I)$ for simplicity in this case.

We have a characterization of
the admissible space of inversive distance circle packings on a weighted triangle.
	
	\begin{lemma}[\cite{Xu MRL}]\label{basic property I of IDCP}
	Let $\triangle v_1v_2v_3$ be a face in $(S, \mathcal{T})$ with three weights $I_1,I_2,I_3\in (1, +\infty)$ defined on edges opposite to the vertices $v_1, v_2, v_3$ respectively.
Let $u: \{v_1, v_2, v_3\}\rightarrow \mathbb{R}$ be a function defined on the vertices, inducing edge lengths by
		\begin{equation}
		\label{length2}
			l_{ij}= \sqrt{e^{2u_i}+e^{2u_j}+2e^{u_i+u_j}I_{k}} =  \sqrt{r_i^2+r_j^2+2r_ir_jI_{k}},
		\end{equation}
where $r_i=e^{u_i}$, $\{i,j,k\}=\{1,2,3\}$.
\begin{enumerate}
  \item[(a)] $l_{12}, l_{13}, l_{23}$ generate a non-degenerate Euclidean triangle if and only if
  \begin{equation}
  \label{definitionQ}
       Q:=\kappa_1^2(1-I^2_{1})+\kappa_2^2(1-I^2_{2})+\kappa_3^2(1-I^2_{3})
		+2\kappa_1\kappa_2\gamma_{3}+2\kappa_1\kappa_3\gamma_{2}+2\kappa_2\kappa_3\gamma_{1}>0,
  \end{equation}
where $\gamma_{i}:=I_{i}+I_{j}I_{k}$ and $\kappa_i:=r_i^{-1}$.
    \item[(b)]  The admissible space $\Omega_{123}$ of inversive distance circle packings $(r_1, r_2, r_3)\in \mathbb{R}^3_{>0}$
on $\triangle v_1v_2v_3$ is
$$\Omega_{123}=\mathbb{R}^3_{>0}\setminus \sqcup_{i=1}^3V_i,$$
where	$\sqcup_{i=1}^3V_i$ is a disjoint union of
		$$V_i=\left\{(r_1, r_2, r_3)\in \mathbb{R}^3_{>0}|\kappa_i\geq \frac{-B_i+\sqrt{\Delta_i}}{2A_i}\right\}$$
		with
		\begin{equation}
		\label{discriminant}
			\begin{aligned}
				A_i=&I^2_{i}-1,\\
				B_i=&-2(\kappa_j\gamma_{k}+\kappa_{k}\gamma_j),\\
				\Delta_i
				=&4(I_1^2+I_2^2+I_{3}^2+2I_1I_2I_{3}-1)(\kappa_j^2+\kappa_{k}^2+2\kappa_j\kappa_{k}I_i).
			\end{aligned}
		\end{equation}

\end{enumerate}
	\end{lemma}



	\subsection{Weighted Delaunay triangulations}\label{subsection weight Delaunay}
	Weighted Delaunay triangulations are natural generalizations of the classical Delaunay triangulations, where the sites generating the corresponding Voronoi decomposition are disks instead of points. It has wide applications in computational geometry.
See \cite{AK, Ed} and others.

\begin{figure}[h]
\begin{overpic}[scale=1]{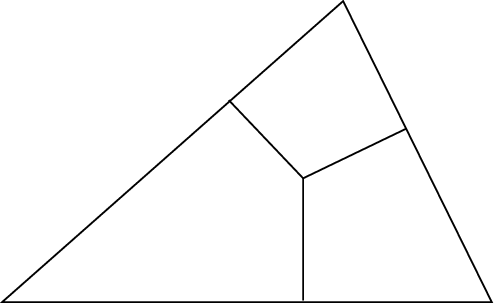}
  \put(-3,3){$v_1$}
    \put(100,3){$v_2$}
        \put(52,2){$C_{12}$}
            \put(33,2){$d_{12}$}
             \put(73,2){$d_{21}$}
             \put(20,26){$d_{13}$}
                \put(50,25){$C_{123}$}
                \put(50,52){$d_{31}$}
                \put(52,40){$h_{13,2}$}
                \put(65,15){$h_{12,3}$}
                \put(66,35){$h_{23,1}$}
                 \put(84,36){$C_{23}$}
                 \put(72,61){$v_{3}$}
                 \put(76,50){$d_{32}$}
                 \put(84,16){$d_{23}$}

                \put(40,44){$C_{13}$}
\end{overpic}
\caption{Sign distances of the geometric center. }
\label{figure1}

\end{figure}

	Assume $r: V\rightarrow (0, +\infty)$ is a non-degenerate inversive distance circle packing on
a weighted triangulated surface $(S, \mathcal{T}, I)$.
Let $\triangle v_1v_2v_3$ be a Euclidean triangle in the plane isometric to a face in $(S, \mathcal{T}, I, r)$.
Then there exists a unique geometric center $C_{123}$ such that its power distances to $v_i$, defined by $|C_{123}-v_i|^2 - r_i^2$, are equal for $i = 1,2,3.$
Projections of the geometric center $C_{123}$ to the lines $v_1v_2, v_1v_3, v_2v_3$ give rise to the geometric centers of these edges, which are denoted by $C_{12}, C_{13}, C_{23}$ respectively.
Please refer to Figure \ref{figure1}.
One can refer to \cite{Glickenstein DCG, Glickenstein JDG, Glickenstein preprint, GT} for more information on the geometric center generated
by discrete conformal structures on manifolds.

Denote $d_{ij}$ as the signed distance of $C_{ij}$ to the vertex $v_i$ and $h_{ij,k}$
as the signed distance of $C_{123}$ to the edge $v_iv_j$.
Glickenstein \cite{Glickenstein JDG} obtained the following identities
	\begin{equation}\label{d}
		d_{ij}=\frac{r_i^2+r_ir_jI_{ij}}{l_{ij}}, \quad h_{ij,k} = \frac{d_{ik} - d_{ij}\cos \theta_i}{\sin \theta_i}.
	\end{equation}
Note that $d_{ij}\in \mathbb{R}_{>0}$ could be defined by (\ref{d})  independent of the existence of the geometric center $C_{ijk}$, and
$h_{ij,k}$ is symmetric in the indices $i$ and $j$, while $d_{ij}$ is not.

For a weighted triangulated surface with a non-degenerate inversive distance circle packing $(S, \mathcal{T}, I, r)$, an interior edge $v_iv_j$ is \textit{weighted Delaunay} if
$h_{ij,k} + h_{ij,l}\geq 0,$
where $\triangle v_iv_jv_k$ and $\triangle v_iv_jv_l$ are two triangles in $F$ sharing the common edge $v_iv_j$.
And $(S, \mathcal{T}, I, r)$ is weighted Delaunay if all the interior edges are weighted Delaunay.

\subsection{Variations of combinatorial curvatures}
The following lemma describes the change of inner angles along PL metrics generated by smooth families of labels on $(S, \mathcal{T}, I)$.

\begin{lemma}[\cite{Guo, Xu AIM}]
\label{derivative angles}
Let $\triangle v_1v_2v_3$ be a face in $(S, \mathcal{T}, I)$ given by Lemma \ref{basic property I of IDCP}.
Suppose that the label $u\in \mathbb{R}^3$ induces a non-degenerate Euclidean triangle $\triangle v_1v_2v_3$, then
\begin{equation}
\label{angle deform}
	    \frac{\partial \theta_i}{\partial u_j}=\frac{\partial \theta_j}{\partial u_i}=\frac{h_{ij,k}}{l_{ij}}, \ \ \
\frac{\partial \theta_i}{\partial u_i}=-\frac{\partial \theta_i}{\partial u_j}-\frac{\partial \theta_i}{\partial u_k}<0,
\end{equation}
		where
		\begin{equation}\label{h_ij,k}
			\begin{aligned}
				h_{ij,k}
			=\frac{r_1^2r_2^2r_3^2}{A_{123}l_{ij}}[\kappa_k^2(1-I_k^2)+\kappa_j\kappa_k\gamma_{i}+\kappa_i\kappa_k\gamma_{j}]
				=\frac{r_1^2r_2^2r_3^2}{A_{123}l_{ij}}\kappa_kh_k
			\end{aligned}
		\end{equation}
		with
$A_{123}=l_{12}l_{13}\sin\theta_1$ and
		\begin{equation}
		\label{h_i}
			\begin{aligned}
h_k=\kappa_k(1-I_k^2)+\kappa_i\gamma_{j}+\kappa_j\gamma_{i}.
			\end{aligned}
		\end{equation}
\end{lemma}


Note that $h_{ij,k}$ given by the equation (\ref{h_ij,k}) is exactly the signed distance of the geometric center $C_{ijk}$ of the triangle $\triangle v_iv_jv_k$ to the edge $v_iv_j$ defined in Subsection \ref{subsection weight Delaunay}.

For a non-degenerate inversive distance circle packing metric $l$ on $(S, \mathcal{T}, I)$,
set $\eta_{ij}^k = h_{ij,k}/l_{ij}$ and  define the \textit{conductance} $\eta: E\rightarrow \mathbb{R}$ for $(S, \mathcal{T}, I, l)$ by
\begin{equation}
\label{definitioneta}
    \eta_{ij} = \begin{cases}
       \eta_{ij}^k + \eta_{ij}^m , & v_iv_j \text{ is an interior edge contained in } \triangle v_iv_jv_k  \text{ and } \triangle v_iv_jv_m;\\
      \eta_{ij}^k , & v_iv_j \text{ is a boundary edge contained in } \triangle v_iv_jv_k.
      \end{cases}
\end{equation}
As a direct corollary of formula (\ref{angle deform}), we have the following variation of combinatorial curvatures.

\begin{corollary}[\cite{Guo, Xu AIM, Xu MRL}]
 Suppose $w(t)*l$ is a family of inversive distance circle packing metrics on $ ({S},\mathcal{T}, I) $ induced by a smooth family of discrete conformal factor $w(t) \in \mathbb{R}^V$. Let $K(t)$ and $\eta(t)$ be the combinatorial curvature and the conductance of $ ({S},\mathcal{T}, I, w(t)*l)$. Then
 \begin{equation}
     \label{curvature}
     \frac{dK_i(t)}{dt} = \sum_{j\sim i}\eta_{ij}(t)\left(\frac{dw_i}{dt} - \frac{dw_j}{dt}\right).
 \end{equation}

\end{corollary}

\subsection{A ring lemma}
Let $P_n$ be a star-shaped $n$-sided polygon in the plane with boundary vertices $v_1,\cdots,v_n$ ordered cyclically ($v_{n+i} = v_i$).
Assume $v_0$ is an interior point of $P_n$ and it induces a triangulation $\mathcal{T}$ of $P_n$ with triangles $\triangle v_0v_iv_{i+1}$.
Then an assignment of radii $r:V(\mathcal{T}) \to \mathbb{R}_{>0}$ is a vector in $\mathbb{R}^{n+1}$.
    \begin{figure}[h]
\begin{overpic}[scale=0.7]{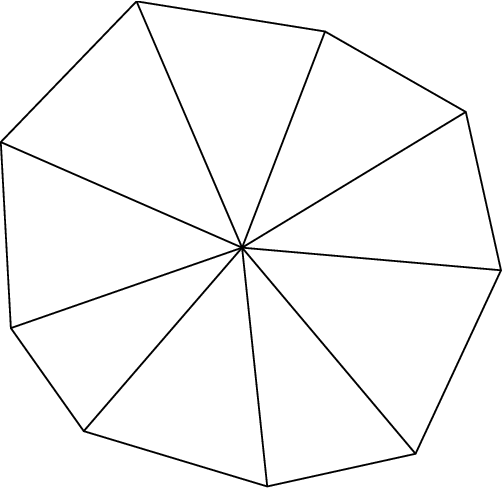}
 \put(90,80){$v_{j+1}$}
  \put(102,45){$v_j$}
  \put(88,10){$v_{j-1}$}
  \put(55,43){$v_0$}
\end{overpic}
\caption{A star triangulation of a polygon.}
   \label{figure3}
\end{figure}

\begin{lemma}\label{ring lemma}
  Let $\mathcal{T}$ be a star triangulation of an $n$-sided polygon $P_n$ with boundary vertices $v_1,\dots, v_n$ and a unique interior vertex $v_0$.
Let  $I:E\rightarrow (1, +\infty)$ be a weight and $r$ be a flat inversive distance circle packing on $(P_n,\mathcal{T}, I)$.
Then there exists a constant $C = C(I, n)>0$ such that $r_0\leq C r_i$ for all $i\in \{1,2, \cdots, n\}$.
\end{lemma}
\begin{proof}
Without loss of generality, we assume $r_0= 1$, otherwise applying a scaling to
the weighed triangulated polygon $(P_n,\mathcal{T}, I, r)$.
Then we just need to prove $Cr_i\geq 1$ for some $C\in \mathbb{R}_{>0}$, which is equivalent to $\kappa_i\leq C$, for all $i\in \{1,2, \cdots, n\}$.

We prove Lemma \ref{ring lemma} by contradiction.
If the result in Lemma \ref{ring lemma} is not true, then there exists a sequence of
inversive distance circle packings $\{r^{(m)}\}_{m=1}^\infty$ with $r^{(m)}_0=1$ on $(P_n,\mathcal{T}, I)$ such that $\lim_{m\rightarrow \infty}\kappa^{(m)}_{i}=+\infty$
for some $i\in \{1,2, \cdots, n\}$.
Without loss of generality, we can assume $i=1$. For the triangle $\triangle v_0v_1v_2$, we set $I_0=I_{12}, I_1=I_{02}$ and $I_2=I_{01}$ for simplicity. As $r^{(m)}$ is a non-degenerate inversive distance circle packing on $(P_n,\mathcal{T}, I)$, we have
\begin{align*}
  (I_2^2-1)(\kappa^{(m)}_2)^2+(I_1^2-1)(\kappa^{(m)}_1)^2+(I_0^2-1)
  -2\kappa^{(m)}_1\gamma_2-2\kappa^{(m)}_2\gamma_1-2\kappa^{(m)}_1\kappa^{(m)}_2\gamma_0< 0
\end{align*}
by Lemma \ref{basic property I of IDCP} (a), which implies
\begin{equation}\label{kappa2 tend infty}
\begin{aligned}
\kappa^{(m)}_2
>& \frac{1}{I_2^2-1}[\kappa^{(m)}_1\gamma_0+\gamma_1
-\sqrt{(I_{0}^2+I_1^2+I_2^2+2I_{0}I_1I_2-1)((\kappa^{(m)}_1)^2+2I_2\kappa^{(m)}_1+1)}]\\
=&\frac{(I_1^2-1)(\kappa^{(m)}_1)^2-2I_2\kappa^{(m)}_1+I_0^2-1}{\kappa^{(m)}_1\gamma_0+\gamma_1
+\sqrt{(I_{0}^2+I_1^2+I_2^2+2I_{0}I_1I_2-1)((\kappa^{(m)}_1)^2+2I_2\kappa^{(m)}_1+1)}}.
\end{aligned}
\end{equation}
Note that $\lim_{m\rightarrow \infty}\kappa^{(m)}_{1}=+\infty$.
We have $\lim_{m\rightarrow \infty}\kappa^{(m)}_2=+\infty$ by (\ref{kappa2 tend infty}),
which is equivalent to $\lim_{m\rightarrow \infty}r^{(m)}_2=0$.
Combining $r^{(m)}_0=1$ with $\lim_{m\rightarrow \infty}r^{(m)}_{1}=\lim_{m\rightarrow \infty}r^{(m)}_2=0$, we have $\lim_{m\rightarrow \infty}(\theta^0_{12})^{(m)}\rightarrow 0$, where $\theta^0_{12}$ is the inner angle of the triangle $\triangle v_0v_1v_2$ at $v_0$.
The same arguments applying to the triangles $\triangle v_0v_iv_{i+1}, i=2, 3, \cdots, n$ subsequently give
$\lim_{m\rightarrow \infty}(\theta^0_{i, i+1})^{(m)}\rightarrow 0$ for all $i=1, 2, \cdots, n$, which implies
$\lim_{m\rightarrow \infty}K^{(m)}_0=2\pi$.
This contradicts the assumption that $v_0$ is a flat interior vertex of $(P_n,\mathcal{T}, I, r^{(m)})$.
\end{proof}

\begin{remark}
  There is another ring lemma in \cite{LXZ} for inversive distance circle packings in the hexagonal triangulated plane with the inversive distance $I\in (-\frac{1}{2}, +\infty)$, the proof of which depends on the combinatorial structure of the hexagonal triangulated plane. The Ring lemma here, i.e. Lemma \ref{ring lemma}, has no restriction on the triangulation, while the inversive distance is required to be in $(1, +\infty)$.
\end{remark}

\subsection{Rigidity of infinite inversive distance circle packings in the plane}
The conformality of the limit of discrete conformal maps $f_n$ in Theorem \ref{conv introduction} is a consequence of the rigidity of infinite inversive distance circle packings in the plane, which is also conjectured by Bowers-Stephenson \cite{BS}.
The following result on the rigidity of infinite inversive distance circle packings in the plane was proved in \cite{LXZ}, which confirms Bowers-Stephenson's rigidity conjecture.

\begin{theorem}[\cite{LXZ}]
\label{infrigidity introduction}
Let $(\mathbb{C}, \mathcal{T}_{st}, I)$ be a weighted hexagonal triangulated plane such that the weight $I: E\rightarrow (1, +\infty)$ is translation invariant.
Assume $l$ is a weighted Delaunay inversive distance circle packing metric on $(\mathbb{C}, \mathcal{T}_{st}, I)$
induced by a constant label.
If $(\mathbb{C}, \mathcal{T}_{st}, I, w*l)$ is a weighted Delaunay triangulated surface isometric to an open set in the plane, then $w$ is a constant function.
\end{theorem}
The rigidity of inversive distance circle packings with prescribed combinatorial curvatures on weighted triangulated compact surfaces has been proved in \cite{Guo,Luo GT,Xu AIM,Xu MRL} based on variational principles. Theorem \ref{infrigidity introduction} provides a result on the rigidity of infinite inversive distance circle packings in the non-compact plane.
In the case of Thurston's circle packings, the rigidity of infinite circle packings in the plane has been explored in \cite{H2, RS, Sch}.

\section{A prescribing curvature problem for inversive distance circle packings}\label{section 3}

To prove Theorem \ref{conv introduction}, we need to establish the existence of discrete conformal maps induced by inversive distance circle packings from a flat polyhedral disk. Following Luo-Sun-Wu's approach in \cite{LSW},
we solve a prescribing combinatorial curvature problem for inversive distance circle packings on it.
This corresponds to the Koebe-Andreev-Thurston Theorem for Thurston's circle packings in the case of circle packings.
As a result, we establish the existence of discrete conformal maps from a flat polyhedral disk to an equilateral triangle for inversive distance circle packings. Furthermore, we have good control on the angle deformations of the triangles in the triangulation, which implies that the discrete conformal maps are $K$-quasiconformal.


Let  $(\mathcal P, \T, l)$ be a flat  polyhedral disk with an equilateral triangulation,  in which all triangles are equilateral. Then the length function $l$ is a constant function on $E$. Given an equilateral Euclidean triangle $\triangle$ in the plane, the \it $n$-th standard subdivision \rm of $\triangle$ is the equilateral triangulation of $\triangle$ by $n^2$ equilateral triangles. Applying this subdivision to each triangle in an equilateral triangulation of a flat polyhedral disk $(\mathcal P, \T, l)$, we obtain its \it $n$-th standard subdivision \rm  $(\mathcal P, \T_{(n)}, l_{(n)})$.
Furthermore, if $l$ is an inversive distance circle packing metric induced by a constant label $u$ and a constant weight $I: E\rightarrow (1, +\infty)$, we require that $l_{(n)}$ is also an inversive distance circle packing metric induced by a constant label $u_{(n)}$ and
a constant weight $I_n: E_{(n)}\rightarrow (1, +\infty)$ taking the same value as $I: E\rightarrow (1, +\infty)$.

\begin{theorem}
\label{exist}
Suppose $(\mathcal P, \T,  l)$ is a flat polyhedral disk with an equilateral triangulation $\T$
such that exactly three boundary vertices $p,q,r$ have combinatorial curvature $\frac{2\pi}{3}$,
and the metric $l$ is an inversive distance circle packing metric
induced by a constant label $u$ and a constant weight $I:E\rightarrow (1, +\infty)$.
	Then for sufficiently large $n$, there is a discrete conformal
	factor $w: V_{(n)}\to \R$ for the $n$-th standard subdivision
	$(\mathcal P, \T_{(n)}, I_{(n)}, l_{(n)} )$  such that
	\begin{enumerate}
	    \item $K_i(w*l_{(n)})=0$ for all $v_i \in V_{(n)}-\{p,q,r\}$,
	    \item $K_i(w*l_{(n)})=\frac{2\pi}{3}$ for all $v_i \in \{p,q,r\}$,
	    \item there is a constant $\theta_0 = \theta_0(I)>0$ independent of $n$ such that all inner angles of triangles in $(\T_{(n)}, w*l_{(n)})$ are in
	the interval $[\theta_0, \pi/2+\theta_0]$.
	\end{enumerate}
\end{theorem}

Note that the underlying metric space of $(\mathcal P, \T_{(n)}, I_{(n)},  w*l_{(n)} )$ is an equilateral triangle,
and $(\mathcal P, \T_{(n)}, I_{(n)},  l_{(n)} )$ is weighted Delaunay for each $n$.

This section is devoted to prove Theorem \ref{exist}. To find the discrete conformal factors in Theorem \ref{exist}, we will construct a system of ordinary differential equations to deform the discrete conformal factors via discrete curvatures at vertices. We first consider such a flow on a standard subdivision of an equilateral triangle in Theorem \ref{dcm triangle}, then use the flow to construct the discrete conformal factor on a flat polyhedral disk required in Theorem \ref{exist}.

In the rest of this section, we assume that any initial polyhedral metric $l$ on $(\mathcal{P},\mathcal{T})$ is an inversive distance circle packing metric induced by a constant weight $I>1$ and a constant label, which is weighted Delaunay. For any discrete conformal factor $w$ on $V$, denote the angle at $v_k$ of a triangle $\triangle v_iv_jv_k$ in the metric $w*l$ as $\theta_{ij}^k(w)$. Similarly, the conductance of $w*l$ defined by the formula (\ref{definitioneta}) is denoted as $\eta(w)$, and the curvature of $w*l$ is denoted as $K(w)$. The notation $a = O(b)$ refers to the fact that $|a| \leq C|b|$ for some constant $C = C(I)>0$.

\subsection{Inversive distance circle packings along flows}
In this subsection, we will solve the following prescribing curvature problem: assume $V_0\subset V$ and the initial curvature of $(\mathcal{P},\mathcal{T}, l)$ is $K^0$. Given a prescribed curvature $K^*$ on $V - V_0$, find a discrete conformal factor $w$ such that $w|_{V_0} = 0$ and $K(w) = K^*$.

Consider a smooth family of discrete conformal factors $w(t)$ satisfying
    \begin{equation}
    \label{flow}
    \left\{ {\begin{array}{lr}
    		{K_i}(w(t)) = (1 - t)K_i^0 + tK_i^*,&v_i \in V - {V_0},\\
    		{w_i}(t) = 0, &v_i \in {V_0},
    \end{array}} \right.\end{equation}
and $w(0) = 0$. This family of $w(t)$, if it exists in the interval $[0,1]$, provides a linear interpolation between the initial curvature $K^0$ and the prescribed curvature $K^*$. Therefore, $w(1)*l$ has curvature $K^*$. Hence, we need to show that this flow exists on $[0,1]$ for some standard subdivisions of $\mathcal{T}$ on $V - {V_0}$. To this end, we recall some basic notions of analysis on graphs.

    Given a graph $(V, E)$,  the set of all oriented edges in $(V,E)$
is denoted  by $\bar{E}$. If $v_iv_j$ is an edge in $E$, we denote it as $i\sim j$. A \it conductance \rm on $G$ is a function $\eta: \bar{E} \to \R_{\geq 0}$ so that $\eta_{ij}=\eta_{ji}$. The following definitions and results are well-known. See \cite{Lov} for details.

\begin{definition} \label{definition lap} Given a finite graph $(V,E)$ with a conductance $\eta$,
	the   gradient  $\triangledown: \R^V \to \R^{\bar{E}}$ is defined by $$ (\nabla f)_{ij} =\eta_{ij}(f_i-f_j),$$
	the  Laplace operator  associated to $\eta$ is the linear map
	$\Delta: \R^V \to \R^V$ defined by
	$$ (\Delta f)_i =\sum_{j \sim i}
	\eta_{ij}(f_i-f_j).$$
\end{definition}

	Given a set $V_0\subset V$ and a function $g: V_0\to \mathbb{R}$, the solution to the Dirichlet problem is a function $f:V\to\mathbb{R}$ satisfying
	$$(\Delta f)_i = 0, \forall v_i\in V-V_0, \text{ and } f|_{V_0} = g.$$

\begin{proposition}\label{mpgraph} Suppose $(V,E)$ is a finite
	connected graph with a conductance $\eta(e)>0$ for any edge $e\in E$. Given a nonempty $V_0 \subset V$ and $g:V_0\to \mathbb{R}$, the solution $f$ to the Dirichlet problem exists. Moreover,
	
	(a) (Maximum principle)
	$ \max_{v_i \in V}
	f_i =\max_{ v_i \in V_0} f_i.$
	
	(b) (Strong maximum principle) If $V-V_0$ is connected and
	$\max_{v_i\in V-V_0} f_i =\max_{v_i \in V_0} f_i$,
	then $f|_{V-V_0}$ is a constant function.
\end{proposition}

Recall that the formula (\ref{definitioneta}) defines a conductance $\eta$ for any inversive distance circle packing on $(S, \mathcal{T}, I)$. If it is weighted Delaunay, then $\eta_{ij}\geq 0$. In the rest of this paper, we assume that the Laplace operator $\Delta$ is induced from this conductance $\eta$ for an inversive distance circle packing.

By the variation formula of curvatures in (\ref{curvature}), we have the following system of ODEs by taking derivative of the equation (\ref{flow}) with respect to $t$
    \begin{equation}
    \label{ode} \left\{ {\begin{array}{lr}
			{{(\Delta w')}_i} = \sum\limits_{j \sim i} {\eta _{ij}}({w'_i} - {w'_j}) = K_i^* - K_i^0,& v_i \in V - {V_0},\\
			{w'_i}(t) = 0, & v_i \in {V_0},
	\end{array}} \right. \end{equation}
with the initial value $w(0) = 0$ and $w'_i = \frac{dw_i}{dt}$. We will show that the solution to the system (\ref{ode}) exists for all $t \in [0,1]$ if $(\mathcal{P},\mathcal{T}, l)$ is chosen carefully. Prior to the existence, we first characterize the maximal interval for the existence of the solution to (\ref{ode}).

Given a weighted triangulated surface $(S, \mathcal{T}, I)$ with an inversive distance circle packing metric $l$, consider the set of discrete conformal factors $W\subset \mathbb{R}^V$ defined by
\begin{equation}\label{space W}
\begin{aligned}
  W = \{w\in \mathbb{R}^V | w*l &\text{ is an inversive distance circle packing metric}\\
   &\text{on $(S, \mathcal{T}, I)$ such that $\eta_{ij}>0$ for all edges}\}.
   \end{aligned}
\end{equation}

\begin{lemma} \label{odeexist} Let $(\mathcal{P}, \mathcal{T}, I)$ be a weighted triangulated surface with an inversive distance circle packing metric $l$ generated by a label $u$. The initial valued problem (\ref{ode})
	defined on $W$ has a unique solution in a maximum interval $[0,	t_0)$ with $t_0>0$ if  $V_0 \neq
	\emptyset$ and $0\in W$. Moreover, if $t_0<\infty$, then either
	$\liminf_{t \to t_0^-} \theta^i_{jk}(w(t)) = 0$ for some angle
	$\theta^i_{jk}$ or $\liminf_{t \to t_0^-} \eta_{ij}(w(t))=0$ for
	some edge $v_iv_j$.
\end{lemma}
\begin{proof}
    The ODE system (\ref{ode}) can be written as
	\begin{equation*} \left\{ {\begin{array}{*{20}{l}}
			{A(w) \cdot w'(t) = b, }\\
			{w(0) = 0,}
	\end{array}} \right. \end{equation*}
    where ${A(w)}$ is a square matrix valued smooth function of $w$, $b$ is a column vector determined by the combinatorial curvature, and ${w'(t)}$ is a column vector. Then ${A(w)}$ is an invertible matrix for a fixed $w\in W$. Indeed, consider the following system of linear equations for a fixed $w$
    \begin{equation} \label{a*f=0} {A(w) \cdot f = 0.}
    \end{equation}
    From (\ref{ode}) we know that equation (\ref{a*f=0}) is equivalent to
    \begin{equation*}   \left\{ {\begin{array}{lr}
    			{{(\Delta f)}_i} = 0,&v_i \in V - {V_0},\\
    			{f_i} = 0,&v_i \in {V_0},
    	\end{array}} \right.
    \end{equation*}
    where ${\eta _{ij}>0}$ for all edges since $w \in W$. The maximal principle in Proposition \ref{mpgraph} implies that ${f = 0}$. Therefore, ${A(w)}$ is invertible.

    As a result, (\ref{ode}) can be written as $w'(t)=A(w)^{-1}b$. Picard's existence theorem for solutions to the ODE systems implies that there exists an interval $[0,{t_0})$ on which  (\ref{ode}) has a solution.

    If $t_0<\infty$ and $t  \nearrow  t_0$, then $w(t)$ leaves every compact set in $W$. Consider subsets $W_{\delta}=\{ w \in W |  \theta^{i}_{jk} \geq \delta,  |w_{i}| \leq \frac{1}{\delta},\eta_{ij} \geq \delta \}$. It is straightforward to check that $W_\delta$ is compact.  Since $w(t)$ leaves every $W_{\delta}$ for each $\delta >0$, one of the following three cases occurs:
    \begin{enumerate}
        \item     $\liminf_{t \to t_0^-} \theta^i_{jk}(w(t))=0$ for some $\theta^i_{jk}$, or
        \item     $\liminf_{t \to t_0^-} \eta_{ij}(w(t))=0$ for some edge , or
        \item     $\limsup_{t \to t_0^-} |w_i(t)|=+\infty$ for some $v_i \in V$.
    \end{enumerate}

    We claim that the case (3) implies the case (1).
    Otherwise, there exists $\delta>0$ such that $\liminf_{t \to t_0^-} \theta^i_{jk}(w(t))>\delta$ for all $\theta^i_{jk}$. Since $w_i'(t) = 0$ for $v_i \in V_0$ along the flow (\ref{ode}), the radius $r_i = e^{w_i+u_i}$ does not change along the flow. Then for any triangle $\triangle v_iv_jv_k$
    with $v_i$ as a vertex, the sine law implies that
    $$\frac{{l_{ij}^2}}{{l_{ik}^2}} \le \frac{1}{{{{\sin }^2}\delta }}, \quad \frac{{l_{ik}^2}}{{l_{ij}^2}} \le \frac{1}{{{{\sin }^2}\delta }},$$
which further implies that ${r_j} \le \frac{{\sqrt I }}{{\sin \delta }}({r_i} + {r_k})$ and ${r_k} \le \frac{{\sqrt I }}{{\sin \delta }}({r_i} + {r_j}) $ by $I>1$ and (\ref{length1 introduction}). Therefore, ${r_j}$ and ${r_k}$ are of the same order.
Specially, $r_j\rightarrow +\infty$ if and only if $r_k\rightarrow +\infty$.
If $w_k$ and $w_j$ go to infinity,  then
$$\cos \theta _{jk}^i = \frac{{l_{ij}^2 + l_{ik}^2 - l_{jk}^2}}{{2{l_{ij}}{l_{ik}}}} = \frac{{r_i^2 + {r_i}{r_j}I + {r_i}{r_k}I - {r_j}{r_k}I}}{{{l_{ij}}{l_{ik}}}} \to -I<-1,$$
   which is impossible. Since the $1$-skeleton of $\mathcal{T}$ is a finite connected graph, we can show inductively that $w_i$ is bounded for all $v_i\in V$, which contradicts the assumption in case (3). This completes the proof for the claim.
\end{proof}

\subsection{Standard subdivisions of an equilateral triangle}
In this subsection, we consider the ODE system (\ref{ode}) when the polyhedral surface is an equilateral triangulation of an equilateral triangle. We will prescribe special curvatures at the boundary vertices such that the discrete conformal maps approximate the power functions in complex analysis. To apply the estimates in network theory, we need to bound the conductance of a weighted Delaunay triangulation as follows.

\begin{lemma}
\label{ratio}
	Let $\triangle v_1v_2v_3$ be a weighted triangle generated by an inversive distance  circle packing $(r_1, r_2, r_3)$ and the weight $I>1$ is a constant. There exists a constant $\theta_0 = \theta_0(I)\in (0, \frac{\pi}{6})$ such that if the three inner angles of the triangle are bounded in $[\pi /6 - \theta_0,  \pi/2 + \theta_0]$, then
	\begin{enumerate}
		\item[(a)] $r_j/r_i \leq 20$ for any two radii $r_i$ and $r_j$,
		\item[(b)] $C\leq \eta_{ij}^k \leq M$ for some constants $C = C(I)>0$  and $M = M(I)>0$.
	\end{enumerate}
\end{lemma}

\begin{proof}
    Set $$\theta_{0} = \min \{\frac{\pi}{1000}, \arcsin \frac{1}{10(20+I)}\}.$$

	To prove part (a), by the angle bound and the sine law,
\begin{equation}\label{lij bdd}
  l_{ij}/l_{ik}\leq 1 /\sin (\pi /6 - \pi /1000)< \sqrt{5}
\end{equation}
 for any two edges in the triangle.
 Without loss of generality, assume that $r_i = 1$. We will prove that $r_j\leq 20$ by contradiction in the following two cases.
	
	If $r_j>20$ and $r_k/r_j\leq 1/5$, then
	$$\frac{l_{ij}^2}{l_{ik}^2} \geq \frac{r_j^2/5 + 2Ir_j+ 1+  4r_j^2/5}{(r_j/5)^2 + 2Ir_j/5+ 1} \geq 5,$$
	which contradicts (\ref{lij bdd}).
	
	If $r_j>20$ and $r_k/r_j> 1/5$, then $r_k>4$ and the inner angle $\theta_{jk}^i$ at $v_i$ is the largest inner angle in $\triangle v_iv_jv_k$.
Note that in this case, we have $I({r_k} + {r_j} - {r_k}{r_j}) + 1<0$, $l_{ij}<\sqrt{I}(r_j+1)$ and $l_{ik}<\sqrt{I}(r_k+1)$.
As a result, by the cosine law, we have
	    $$
	\begin{aligned}
		\cos \theta _{jk}^i &  = \frac{{l_{ij}^2 + l_{ik}^2 - l_{jk}^2}}{{2{l_{ij}}{l_{ik}}}} \\
&= \frac{{I({r_k} + {r_j} - {r_k}{r_j}) + 1}}{{{l_{ij}}{l_{ik}}}} \\
		&  < \frac{{ - I({r_k} + {r_j} + {r_k}{r_j} + 1) + 2I({r_k} + {r_j}) + I + 1}}{{I({r_k} + {r_j} + {r_k}{r_j} + 1)}}
		\\
		& \le   - 1 + \frac{{2I({r_k} + {r_j}) + 2I}}{{I({r_k} + {r_j} + {r_k}{r_j} + 1)}}\\
& \le \frac{{ - 11}}{{21}}.\\
	\end{aligned}
	$$
	This contradicts that the angle bound is $[\pi /6 - \theta_0,  \pi/2 + \theta_0]$.

	
	To prove part (b), the definition of $\eta_{ij}^k$ and the formula (\ref{h_ij,k}) implies
	$\eta_{ij}^k = \frac{h_{ij,k}}{l_{ij}} = \frac{r_i^2r_j^2r_kh_k}{Al^2_{ij}},$
	where $A=l_{ij}l_{ik}\sin \theta^i_{jk}$. The sign of $\eta_{ij}^k$ is determined by $h_k$. We will show
	\begin{equation}
	\label{radius}
	    	r_kh_k = (1+I)(1 + I(\frac{r_k}{r_i} +\frac{r_k}{r_j} - 1))\geq \frac{1+I}{4}>0.
	\end{equation}
	We just need to check the case that $r_k/r_i\leq 1$ and $r_k/r_j\leq 1$.  If both $r_k/r_i\geq 1/2$ and $r_k/r_j\geq 1/2$, then $r_kh_k\geq 1+I$. Hence, we only need to consider the situation that  $r_k/r_i<1/2$ or $r_k/r_j<1/2$. By the angle bound and cosine law, we have
	$$-\frac{1}{5(20+I)}l_{ij}^2 \leq l_{jk}^2 + l_{ik}^2 -l_{ij}^2 .$$
	This is equivalent to
		$$I(r_ir_k + r_jr_k - r_ir_j ) \geq - r_k^2 - \frac{1}{10(20+I)}(r_i^2 + r_j^2 + 2Ir_ir_j),$$
which implies
	 $$1 + I(\frac{r_k}{r_i} +\frac{r_k}{r_j} - 1) \geq 1 - \frac{r_k^2}{r_ir_j} - \frac{1}{10(20+I)}(\frac{r_i}{r_j} + \frac{r_j}{r_i} + 2I) \geq \frac{4}{5} - \frac{r_k^2}{r_ir_j},$$
where the results in part (a) of Lemma \ref{ratio} is used in the last inequality.
	Then by the formula (\ref{h_i}) of $h_k$, we have
	$$r_kh_k = (1+I)(1 + I(\frac{r_k}{r_i} +\frac{r_k}{r_j} - 1)) \geq (1+I)(\frac{4}{5} - \frac{r_k^2}{r_ir_j}).$$
	 Therefore, under the assumption that $r_k/r_i\leq 1$ and $r_k/r_j\leq 1$, we have $r_kh_k\geq 3(1+I)/10> (1+I)/4$ when  $r_k/r_i<1/2$ or $r_k/r_j<1/2$.

	The sine law implies that $l_{ij}^2/50\leq A\leq 5l_{ij}^2$. Combining with part (a) of Lemma \ref{ratio}, we can find two constants $M=M(I)$ and $C=C(I)$ such that
	$$M(I) \geq \frac{r_i^2r_j^2 (1+I)(1+40I)}{Al^2_{ij}}\geq \frac{r_i^2r_j^2r_kh_k}{Al^2_{ij}}\geq \frac{r_i^2r_j^2 (1+I)}{4Al^2_{ij}} \geq C(I) > 0.$$ \end{proof}

\begin{theorem}
\label{dcm triangle} Let $\mathcal{P}=\triangle ABC$ be an
	equilateral triangle,  $\T_{(n)}$ be the $n$-th standard subdivision of
	$\mathcal{P}$, $l$ be an inversive distance circle packing metric on $(\mathcal{P}, \mathcal{T}_{(n)})$ induced by
 a constant weight $I$ and a constant label.  Set
	$$V_0=\{ v \in V | v \text{ is in the edge  $BC$ of the triangle }\Delta ABC\}.$$
	Given any $\alpha \in [\frac{\pi}{6},
	\frac{\pi}{2}]$, there exists a smooth family of discrete conformal factors $w(t) \in
	\R^V$ for $t\in [0,1]$ such that $w(0)=0$ and  $w(t)* l $ is an inversive distance circle packing metric
	 on $\T_{(n)}$ with curvature $K(t)=K(w(t)* l)$ satisfying
	\begin{enumerate}
	    \item 	   	 $K_A(t)= -t\alpha+(2+t)\frac{\pi}{3}$,
        \item 	    $K_i(t)=0$  \text{for all $v_i \in V-\{A\}\cup V_0$,}
	
        \item 	    $w_i(t)=0$ \text{for all $v_i \in V_0$,}
	
	\item all inner angles $ \theta^i_{jk}(t)$ in the metric $w(t)* l$
	are in the interval  $$ [\frac{\pi}{3}-|\alpha-\frac{\pi}{3}|,
	\frac{\pi}{3}+|\alpha-\frac{\pi}{3}|] \subset [\frac{\pi}{6},
	\frac{\pi}{2}],$$
	
	
	\item   for $v_i\neq A$,
	$$|K_i(t) -K_i(0)| = O(\frac{1}{\sqrt{ \ln(n)}}).$$Moreover,
	\begin{equation}\label{totalcur} \sum_{v_i \in V_0}
		|K_i(t)-K_i(0)|\leq \frac{\pi}{6}.  \end{equation}
	\end{enumerate}
\end{theorem}

Notice that the angle at the vertex $A$ is $t\alpha+(1-t)\pi/3$ along $w(t)$, and curvatures of vertices stay zero except vertices in $BC$ and the vertex $A$.
The piecewise linear map from $(\mathcal{P}, \mathcal{T}_{(n)}, l)$ to $(\mathcal{P}, \mathcal{T}_{(n)}, w*l)$
determined by Theorem \ref{dcm triangle} is a discrete analogue of the analytic function $f(z) = z^{3\alpha/\pi}$. This construction works for any $n$-th subdivision of equilateral triangulations.

To prove Theorem \ref{dcm triangle}, we need the following two estimates for solutions to the Dirichlet problem when the graph is an equilateral triangulation of a polygonal disk.

\begin{lemma}[\cite{LSW}, Lemma 5.8] \label{estimate1}
	Assume $\Delta ABC,n,\mathcal T,V_0$ are as given in Theorem \ref{dcm triangle}.
	Let $\tau: \T \to \T$ be the involution
	induced by the reflection of $\Delta ABC$ about the angle bisector
	of $\angle BAC$ and $\eta: E \to \R_{\geq 0}$ be a conductance so
	that $\eta \tau =\eta$ and $\eta_{ij}=\eta_{ji}$. Let $\Delta:
	\R^V \to \R^V$ be the Laplace operator defined by $(\Delta f)_i
	=\sum_{ j \sim i} \eta_{ij}(f_i -f_j)$. If $f \in \R^V$ satisfies
	$(\Delta f)_i=0$ for $v_i \in V-\{A\}\cup V_0$ and $f|_{V_0}=0$,
	then for all edges $v_iv_j$, the gradient $(\triangledown
	f)_{ij}=\eta_{ij}(f_i-f_j)$ satisfies
	\begin{equation}\label{90e} |(\triangledown
	f)_{ij}| \leq
		\frac{1}{2}|\Delta(f)_A|.
	\end{equation}
\end{lemma}

\begin{lemma}[\cite{LSW}, Lemma 5.9] \label{estimate2}
	Assume $\Delta ABC,n,\mathcal T,V_0$ are as given in Theorem \ref{dcm triangle}.
	Let $\eta: E(\T) \to [\frac{1}{M}, M]$ be a conductance
	function  for some $M>0$ and  $\Delta$ be the Laplace operator on
	$\R^V$ associated to $\eta$. If $f:V \to \R$ solves the Dirichlet
	problem $(\Delta f)_i=0, \forall v_i \in V-\{A\}\cup V_0$,
	$f|_{V_0}=0$ and $(\Delta f )_A=1$, then for all $u \in V_0$,
	$|(\Delta f)_u|\leq \frac{20M}{\sqrt{\ln n}}.$
\end{lemma}

\begin{proof}[Proof of Theorem \ref{dcm triangle}]
    We will prove Theorem \ref{dcm triangle} by considering the ODE system  (\ref{ode}) for $(\triangle ABC, \mathcal{T}_{(n)}, l)$ when $n$ is sufficiently large. The prescribed curvature $K^*$ is
$$K_A^* = \pi  - \alpha, K_i^* = 0, v_i \in V - {V_0} \cup \{ A\}.$$
	The initial curvature ${K^0}$ is
$$K_i^0 = K_i(0) = \frac{2\pi }{3},v_i \in \{ A,B,C\},
			K_i^0 = K_i(t) = 0,v_i \in V - \{ A,B,C\}.$$
Then the ODE system  (\ref{ode}) could be written as
	\begin{equation} \label{flowtriangle}  \left\{ {\begin{array}{lr}
				{K_A}'(t) = {{(\Delta w')}_A} =\frac{\pi }{3} - \alpha, &\quad \\
				{K_i}'(t) = {{(\Delta w')}_i} = K_i^* - K_i^0 = 0,&v_i \in V - {V_0} \cup \{ A\}, \\
				{w'_i}(t) = 0,& v_i \in {V_0},
		\end{array}} \right.
	 \end{equation}
	 with the initial value ${w_i}(0) = 0$. It is straightforward to check that $w(0)\in W$, where $W$ is defined by (\ref{space W}).

   Then there exists a maximum ${t_0} > 0$ such that a solution $w(t)$ to (\ref{flowtriangle}) exists. Moreover, Lemma \ref{odeexist} implies that there exists a maximal time $s_0$ such that $w(t) \in W$ and the statement (4) holds for $t\in [0, s_0)$. We will prove $s_0\geq 1$. Moreover, $w(1)$ exists and $w(1)*l$ is an inversive distance circle packing metric satisfying (1)-(4) in Theorem \ref{dcm triangle}. Without loss of generality, we assume that $s_0<\infty$.

   \textbf{Claim}: For any inner angle $\theta^i_{jk}$ and  $t\in [0, s_0)$, we have
      \begin{equation}\label{theta bound}
      |\theta _{jk}^i(t) - \frac{\pi }{3}| \le {t}|\alpha  - \frac{\pi }{3}|.
    \end{equation}
   To prove this claim, notice that $\alpha\in [\pi/6, \pi/2]$ implies $\theta^i_{jk}(t)\in [\pi/6, \pi/2]$ for $t\in [0, s_0)$ by the statement (4). By Lemma \ref{ratio}, $\eta_{ij}^k(t)>C(I)>0$ for any triangle $\triangle v_iv_jv_k$. Then
       $$|{(\nabla w')_{ij}}| = {\eta _{ij}}|{w'_i} - {w'_j}| \ge \eta _{ij}^k|{w'_i} - {w'_j}|.$$
    By Lemma \ref{estimate1} and  ${\frac{{d{K_i}}}{{dt}} = {{(\Delta w')}_i}}$, we obtain
    \begin{equation*}  |{(\nabla w')_{ij}}| \le \frac{1}{2}|{(\Delta w')_A}| = \frac{1}{2}\Big|\frac{{d{K_A}}}{{dt}}\Big| = \frac{1}{2}|\alpha  - \frac{\pi }{3}|.
    \end{equation*}
    By the formula (\ref{angle deform}), we have
    \begin{equation*} \Big|\frac{{d\theta _{jk}^i}}{{dt}}\Big| \le \eta _{ik}^j|{w'_i} - {{w'}_k}| + \eta _{ij}^k|{w'_i} - {w'_j}| \le |{(\nabla w')_{ik}}| + |{(\nabla w')_{ij}}| \le |\alpha  - \frac{\pi }{3}|.
    \end{equation*}
    Then for all $t \in [0,s_0)$,  	
    \begin{equation*}  |\theta _{jk}^i(t) - \frac{\pi }{3}| = |\theta _{jk}^i(t) - \theta _{jk}^i(0)| \le \int_0^t \Big| \frac{{d\theta _{jk}^i(t)}}{{dt}}\Big|dt \le t|\alpha  - \frac{\pi }{3}|.
    \end{equation*}

   Now it is not hard to show that $s_0\geq 1$ from the claim. Notice that $\liminf_{t\to s_0^-} \theta^i_{jk}(w(t)) \geq \pi/6$ and $ \liminf_{t \to s_0^-} \eta_{ij}(w(t))>0$ for all edges by Lemma \ref{ratio}. Therefore, as $t\to s_0^-$, for some $\theta^i_{jk}$,
\begin{equation}\label{theta limit as t tends s0}
  \limsup_{ t \to s_0^-}
    		|\theta^i_{jk}(w(t))-\frac{\pi}{3}|= |\alpha-\frac{\pi}{3}|
\end{equation}
    by the definition of $s_0$ and Lemma \ref{odeexist}.
   If $s_0<1$, then for all the inner angles, we have
   \begin{equation*}  \limsup_{ t \to s_0^-}|\theta _{jk}^i(t) - \frac{\pi }{3}| \le s_0|\alpha  - \frac{\pi }{3}| < |\alpha  - \frac{\pi }{3}|
    \end{equation*}
    by (\ref{theta bound}), which contradicts (\ref{theta limit as t tends s0}). Therefore, $s_0\geq 1$.

   Notice that $\theta_{jk}^i(t) \in [\pi/6, \pi/2]$ for all inner angles and $\eta_{ij}(t)\geq C(I)>0$ by Lemma \ref{ratio} when $t\in [0, s_0)$. This means that $w(t)*l$ is non-degenerate and strictly weighted Delaunay for $t\in [0, s_0)$. Therefore, Lemma \ref{odeexist} implies that $1\leq s_0<t_0$. Then $w(1) \in W$. By continuity, the metric $w(1)*l$ is a non-degenerated inversive distance circle packing
   metric satisfies (1)-(4) in Theorem \ref{dcm triangle}.

	Finally, we use Lemma \ref{estimate2} to prove the last statement in Theorem \ref{dcm triangle}.
Notice that by Lemma \ref{ratio}, $0< C\leq \eta_{ij} \leq M$, where $C=C(I), M=M(I)$. Applying Lemma \ref{estimate2} to the function $f = \frac{{dw(t)}}{{dt}}/|\alpha  - \pi /3|$, we obtain
	\begin{equation*}\Big|\frac{{d{K_i}(t)}}{{dt}}\Big| = |{(\Delta w')_i}| = |\alpha  - \pi /3| \cdot |{(\Delta f)_i}| = O(\frac{1}{{\sqrt {\ln (n)} }}),v_i \in V_0 .
	\end{equation*}	
    Then for $v_i \ne A$, we have for $t\in[0, 1]$,
    $$|{K_i}(t) - {K_i}(0)| \le \int_0^t \Big| \frac{{d{K_i}(t)}}{{dt}}\Big|dt = O(\frac{1}{{\sqrt {\ln (n)} }}).$$
    Moreover, if $\alpha = \pi/3$, then (\ref{totalcur}) is automatically true since the flow would be a constant flow
    by Proposition \ref{mpgraph}. Hence we assume $\alpha\neq \pi/3$. We claim that $ {{w'}_A}(t) \ne 0$ for $t \in [0,{t_0})$. Otherwise, ${{w'}_A}(s) = 0$ for some $s \in [0,{t_0})$. Applying the maximum principle, i.e. Proposition \ref{mpgraph}, to the following Dirichlet problem
	$${\left\{ {\begin{array}{lr}
				{{(\Delta w'(s))}_i= 0}, &v_i \in V - \{ A\}  \cup {V_0},\\
				{w'_i}(s) = 0, &v_i \in \{ A\}  \cup {V_0},
		\end{array}} \right.}$$
	 we obtain ${w'_i}(s) = 0$ for all $v_i \in V$, which implies ${(\Delta w')_A} (s)= 0$. This contradicts ${(\Delta w')_A} (s)=\alpha  - \frac{\pi }{3} \neq 0$. This completes the proof of the claim.
Furthermore, applying the maximal principle, i.e. Proposition \ref{mpgraph}, to (\ref{flowtriangle}) again shows that $w'_A(t)$ and $w'_i(t)$ have the same sign. Note that for $v_i\in V_0$, $w_i'(t) = 0$. We have
	 $$w'_A(t)K'_i(t) =w'_A(t) \sum_{i\sim j}\eta_{ij}(w_i' - w_j') =   -\sum_{i\sim j}\eta_{ij}w'_j(t)w'_A(t) \leq 0, $$
	 which implies $({K_i}(t) - {K_i}(0)){{w'}_A}(t) \le 0$. By the discrete Gauss-Bonnet formula (\ref{discrete Gauss-Bonnet formula}), we have
    $${K_A}(t) + \sum\limits_{v_i \in {V_0}} {{K_i}} (t) = {K_A}(0) + \sum\limits_{v_i \in V} {{K_i}} (0) = 2\pi . $$
	Since ${K_i}(t) - {K_i}(0)$ have the same sign for all $v_i \in {V_0}$, we conclude that for $t\in[0,1]$,
    $$\sum\limits_{v_i \in {V_0}} | {K_i}(t) - {K_i}(0)| =|\sum\limits_{v_i \in {V_0}}  ({K_i}(t) - {K_i}(0))| = |{K_A}(t) - {K_A}(0)| = |t(\alpha  - \frac{\pi }{3})| \le \frac{\pi }{6}.$$
\end{proof}	

\subsection{Proof of Theorem \ref{exist}}
    There are two steps to find the discrete conformal factor required in Theorem \ref{exist}. In the first step, we construct a discrete conformal factor by Theorem \ref{dcm triangle} where the triangles contain boundary vertices of nonzero curvatures. This step will diffuse the curvature of boundary vertices of polyhedral disk $\mathcal{P}$ to interior vertices such that curvatures are small if the subdivision is sufficiently dense. In the second step, we eliminate these small curvatures using a flow similar to (\ref{ode}). We need the following lemma in the second step.

  \begin{lemma}[\cite{LSW}, Proposition 5.10] \label{estimate4} Suppose $(\mathcal P, \T', l)$ is
	polygonal disk  with an equilateral triangulation and
	$\T$ is the $n$-th standard subdivision of the triangulation
	$\T'$ with $n \geq e^{10^6}$. Let $\eta: E=E(\T) \to [\frac{1}{M},
	M]$ be a conductance function with $M>0$ and $\Delta: \R^V \to \R^V$ be the
	associated Laplace operator.  Let
	$V_0 \subset V(\T)$ be a thin subset such that for all $v \in V$ and
	$m \leq n/2$, $|B_m(v) \cap V_0| \leq 10m$.   If $f: V \to \R$
	satisfies $(\Delta f)_i=0$ for $v_i \in V-V_0$, $|(\Delta f)_i| \leq
	\frac{M}{\sqrt{\ln(n)}}$ for  $v_i \in V_0$ and $\sum_{v_i \in V_0}
	|(\Delta f)_i| \leq M$, then for all edges $v_jv_k$ in $\T$,
	$$|f_j-f_k| \leq \frac{200M^3}{\sqrt{\ln(\ln(n))}}.$$
    \end{lemma}

\begin{proof}[Proof of Theorem \ref{exist}]
    We call each boundary vertex of $\mathcal{P}$ other than ${p, q, r}$ \emph{corner} if it has nonzero curvature. Denote the set of corners as $V_c$. Notice that by the assumption on $\mathcal{P}$, each vertex in $V_c$ has degree $m = 3, 5$ or $6$. Moreover, the standard subdivision of each triangle of $\mathcal{P}$ does not introduce new corners. Thus, the cardinality $|V_c|$ of $V_c$ is independent of the subdivision $\mathcal{T}_{(n)}$ of $\mathcal{T}$.

    Let  ${B_{[n/3]}}(v)$ be the combinatorial ball in $\mathcal{T}_{(n)}$ centered at $v\in V_c$ with radius $ [n/3]$ where $[x]$ is the integer part of a real number $x$. Notice that these ball are disjoint in $\mathcal{T}_{(n)}$. Each ${B_{[n/3]}}(v)$ consists of $m-1$ copies of $[n/3]$-th subdivision of equilateral triangles $\triangle_1^v, \dots, \triangle_{m-1}^v$.

    \textbf{Step 1}: For every $v\in V_c$, we will deform its curvature to zero. In particular, we apply Theorem \ref{dcm triangle}
    to $\triangle_1^v, \dots, \triangle_{m-1}^v$ with $\alpha = \pi/(m-1)\in [\frac{\pi}{6}, \frac{\pi}{2}]$. It produces a discrete conformal factor $w_i$ on $\triangle_i^v$ for each $i = 1, \dots, m-1$. Notice that the discrete conformal factor on $\mathcal{T}_{(n)}$ in Theorem \ref{dcm triangle} depends only on $\alpha$. Then discrete conformal factor $w_i$ are identical on each $\triangle_i$. By symmetry, we can glue them together to form a discrete conformal factor $\bar w$ on $\mathcal{T}_{(n)}$. Specifically, the value of $\bar w$ on  ${B_{[n/3]}}(v)$ for $v\in V_c$ is determined by Theorem \ref{dcm triangle}, and the values of $\bar w$ on other vertices are zero.

Let $\bar K$ be the curvature of  inversive distance circle packing metric $\bar l = \bar w*l$.
Let $K$ be the curvature of the target equilateral triangle with $K_i=0$ for all $v_i \in V_{(n)}-\{p,q,r\}$ and $K_i=\frac{2\pi}{3}$ for all $v_i \in \{p,q,r\}$. Then Theorem \ref{dcm triangle} implies that
    \begin{enumerate}
        \item    $\bar K_i = K_i$ for all vertices $v_i$ in the set of $V_k := \{ v_i|{d_c}(v_i,v) \ne [n/3],v \in {V_c}\}$,
        \item  $\bar w_i = 0$ for all $v_i \notin { \cup _{v \in {V_c}}}{B_{[n/3]}}(v)$,
    \item   all inner angles at $v \in {V}$ satisfy $\theta _{ij}^v \in [\frac{\pi }{6},\frac{\pi }{2}]$,
    \item for all vertices $ v_i \notin V_k$, $|{{\bar K}_i} - K_i| = O(\frac{1}{{\sqrt {\ln (n)} }})$,
    \item
        $\sum\limits_{v_i \in V} | {{\bar K}_i} - K_i| \le \frac{2\pi N}{3}$, where $N$ denote the number of corners.
    \end{enumerate}
    Notice that the set $V_k$ is the union of the sets $V_0$ given by Theorem \ref{dcm triangle} for each $v\in V_c$. Statement (1) and (2) are immediate by the construction. Statement (3), (4), and (5) are immediate by Theorem \ref{dcm triangle}.

    \textbf{Step 2}: we construct a flow to deform the curvatures of vertices in $V_k$ to be zero when the subdivision is sufficiently dense. Specifically, consider the following ODE system on $\mathcal{T}_{(n)}$
       \begin{equation} \label{odestep2} \left\{ {\begin{array}{lr}
   			\frac{{d{K_i}(w(s)*\bar l)}}{{dt}} = K_i -  {{{\bar K}_i}} ,& v_i \in V - \{ p,q,r\},\\
   			{w_i}(s) = 0,& v_i \in \{ p,q,r\},
   	\end{array}} \right.
   \end{equation}
    with initial value $w(0) = 0$. The idea  is the same as that of (\ref{ode}). Namely, we linearly interpolate the initial curvature $\bar K$ and the target curvature $K$. By Lemma \ref{odeexist}, there exists a maximal $s_0>0$ such that the solution $w(s)$ to (\ref{odestep2}) exists and $w(s)*\bar l$ satisfies  that  on $[0, s_0)$, all inner angles at $v \in {V}$, $\theta _{ij}^v \in [\frac{\pi }{6} - \theta_0,\frac{\pi}{2} + \theta_0]$, where $\theta_0$ is the parameter given by Lemma \ref{ratio}.


    Now we apply Lemma \ref{estimate4} to estimate the angle deformation along the flow (\ref{odestep2}).
    Set $V_B=V_{(n)}\setminus V_k$.
    First notice that $V_B$ is a thin set in $V_{(n)}$ of $\mathcal{T}_{(n)}$. In particular,  $|{B_r}(v_i) \cap {V_B}| \le 10r$ for
    any $v_i\in V_{(n)}$ and any $r \le n/3$. Moreover, by (4) and (5) in Step 1, we obtain
    $$\sum\limits_{v_i \in {V_B}} | {K_{i}'}| =\sum\limits_{v_i \in {V_B}} | {(\Delta w')_i}| \le \sum\limits_{i\in V} | {{\bar K}_i} - K_i| \le \frac{{2\pi N}}{3},$$ and
    $$|{K_{i}'}|= |{(\Delta w')_i}|  \le |{{\bar K}_i} - K_i| = O(\frac{1}{{\sqrt {\ln (n)} }}) ,v_i \in {V_B} .$$
    Lemma \ref{ratio} implies that  $f = w'$ along the flow (\ref{odestep2}) satisfies the conditions in Lemma \ref{estimate4}. Therefore, we obtain that if $i \sim j$, then
    $$ |{w'_i}(s) - {w'_j}(s)| = O(\frac{1}{{\sqrt {\ln(\ln (n)}) }}).$$
As a result,
    $$ \Big|\frac{{d\theta _{ij}^k}}{{ds}}\Big| \le |\eta _{jk}^i({w'_j} - {{w'}_k})| + |\eta _{ik}^j({w'_i} - {{w'}_k})|  = O(\frac{1}{{\sqrt {\ln(\ln (n)}) }}), $$
    where Lemma \ref{ratio} is used in the last equality.
    For all $s \in [0,{s_0})$ and sufficiently large $n$,
    \begin{equation*}|\theta _{ij}^k(w(s)) - \theta _{ij}^k(0)| \le \int_0^s \Big| \frac{{d\theta _{ij}^k(w(s))}}{{ds}}\Big|ds = O(\frac{1}{{\sqrt {\ln(\ln (n)}) }}) \le \theta_0 s_0.
    \end{equation*}
    We claim that ${s_0} > 1$. Otherwise, we can extend the solution to (\ref{odestep2}) to $[0, s_0+\epsilon)$ for some small $\epsilon>0$, which contradicts the maximality of $s_0$. Set $w^*=w(1)$ and $w = \bar w + w^*$. Then the curvature of the inversive distance circle packing metric $w*l$ is $$K(0)+\int^1_0 K'(s)ds=\bar{K}+(K-\bar{K})=K.$$
    This implies that the discrete conformal factor $w = \bar w +w^*$ produces the discrete conformal map from $(\mathcal{P}, \mathcal{T}_{(n)}, l_{(n)})$ to the equilateral triangle.
    \end{proof}

\section{The convergence of inversive distance circle packings}\label{section 4}


To prove Theorem \ref{conv introduction}, we first recall the following three theorems on the extension and convergence of quasiconformal maps.
\begin{theorem}[\cite{Al}, Corollary in Page 30]
\label{extension}
If $f: \D \to \Omega$ is a $K$-quasiconformal map from the open unit disk $\mathbb{D}$ onto a Jordan
domain $\Omega$, then $f$ extends continuously to a
homeomorphism $\overline{f}: \overline{\D} \to
\overline{\Omega}$.
\end{theorem}
The following theorem is a simple consequence of Lemma 2.1 and Theorem 2.2 in \cite{Leh}.
\begin{theorem}
\label{compactness}
If $f_n:\mathbb{D}\to\Omega_n$ is a sequence of $K$-quasiconformal maps such that $\Omega_n$ is uniformly bounded, then every subsequence of $f_n$ contains a subsequence that converge locally uniformly. Moreover, the limit of this subsequence is a $K$-quasiconformal map or a constant map.
\end{theorem}

A sequence of Jordan curves
$J_n$ in $\C$ converge uniformly to a Jordan curve $J$
in $\C$ if there exist homeomorphisms $\phi_n: \mathbb S^1 \to
J_n$ and $\phi: \mathbb S^1 \to J$ such that $\phi_n$ converge
uniformly to $\phi$.

\begin{theorem}[\cite{Pal}, Corollary 1] \label{pal} Assume that $\Omega_n$ is a sequence of Jordan domains such that $\partial \Omega_n$ converge uniformly to  $\partial \Omega$.  If $f_n: \D \to \Omega_n $ is a  $K$-quasiconformal map for each $n$, and the sequence \{$f_n$\}
converge to a $K$-quasiconformal map $f: \D \to \Omega$ uniformly
on compact sets of $\D$, then $\overline{f_n}$ converge to
$\overline{f}$ uniformly on $\overline{\D}$.
\end{theorem}



\begin{proof}[Proof of Theorem \ref{conv introduction}]
By taking the intersection of scalings of the standard hexagonal triangulation in the plane with $\Omega$, we can construct a sequence of nested polygonal disks $\Omega_n$ such that $\partial \Omega_n$ converge uniformly to $\partial \Omega$ and there are three boundary vertices $p_n,q_n, r_n \subset
\partial \Omega_n$ such that $\lim_n p_n =p$, $\lim_n q_n =q$ and
$\lim_n r_n =r$. By adding or subtracting boundary vertices if necessary, we can assume that the curvatures at $p_n,
q_n, r_n\in \partial \Omega_n$ are $\frac{2\pi}{3}$ and the curvatures at all
other boundary vertices of $\Omega_n$ are not $\frac{2\pi}{3}$.

By Theorem \ref{exist}, we produce some standard subdivision $\T_n$ of
$\Omega_n$ and some discrete conformal factors $w_n$ such that $(\Omega_n, \T_n,
w_n*l_{st})$ is isometric to the unit equilateral triangle $(\triangle ABC, \T_n)$, where
$A,B,C$ correspond to $p_n,q_n,r_n$ respectively.
Let $f_n: (\triangle ABC,\mathcal{T}_n, (A,B,C)) \to (\Omega_n, \T_n, (p_n, q_n, r_n))$ be the discrete conformal map induced by the correspondence of triangulations.
Let $\bar{f}$ be  the Riemann mapping from $\triangle ABC$ to $\overline{\Omega}$
 sending $A, B, C$ to $p, q, r$ respectively. We claim that $f_n$ converges uniformly to
$\bar{f}$ on $\triangle ABC$. 

By Theorem \ref{exist},  all angles of triangles in  $(\triangle ABC, \T_n, w*l_{(n)})$ are
at least $\epsilon_0>0$.  Then the discrete conformal maps $f_n$ are
$K$-quasiconformal from $int(\triangle ABC)$ to $int(\Omega_n)$ for some constant $K$ independent of $n$ and continuous from $\triangle ABC$ to $\Omega_n$. Let $\mathring{f}_n$ be the restriction of $f_n$ in $int(\triangle ABC)$. Theorem \ref{compactness} implies that every convergent subsequence of $\{\mathring f_n\}$ converge to a $K$-quasiconformal map $\mathring{g}$ from $int(\triangle ABC)$ to $int(\Omega)$. Since $\Omega = \cup_n\Omega_n$, $\mathring{g}$ is onto $int(\Omega)$. Theorem \ref{extension} implies that $\mathring{g}$ extends to a homeomorphism $g:\triangle ABC \to \Omega$. Theorem \ref{pal} implies that $f_n$ converge uniformly to $g$ on $\triangle ABC$. It is straightforward to check that $g(A) = p$, $g(B) = q$, and $g(C) = r$.

Notice that the Riemann mapping $\bar f$ is the only continuous extension of a \textit{conformal} map from $int(\triangle ABC)$ to $\Omega$ with $\bar{f}(A) = p$, $\bar{f}(B) = q$, and $\bar{f}(C) = r$. This means that if we can show $g$ is conformal, then $g = \bar{f}$ and all limits of convergent subsequences of $\{f_n\}$ are $\bar{f}$. This will complete the proof of $f_n\to \bar{f}$ uniformly on $\triangle ABC$.

The conformality of $g$ follows from Theorem \ref{infrigidity introduction} by the same argument as the Hexagonal Packing Lemma in \cite{RS}. We briefly repeat the arguments here for completion.
For a vertex $v_0\in \mathcal{T}_{st}$, let $B_n$ be the $n$-ring neighborhood of $v_0$ in $\mathcal{T}_{st}$. Then $B_n$ is a finite simplicial complex whose underlying space is a topological disk $\mathbb{D}$. Assume that $l_n$ is a flat inversive distance circle packing on $B_n$ with the constant weight $I$. Let $s_n$ be the maximal ratio of radii of two adjacent circles in $l_n$ of $B_n$. Lemma \ref{ring lemma} implies that $s_n$ is uniformly bounded by some constant $C(I)$. As $n\to \infty$, we can pick a convergent subsequence of $(\mathbb{D}, B_n, I, l_n)$, still indexed with $n$, such that all circles converge geometrically. We claim that $\lim_n{s_n} = 1$. Otherwise, as $n\to \infty$, the limit produces an inversive distance circle packing on $\mathcal{T}_{st}$ such that circles have different sizes. This contradicts  the fact that $w$ is a constant in Theorem \ref{infrigidity introduction}.

As $n\to \infty$, the arguments above show that $s_n$ of $\mathcal{T}_n$ goes to $1$. Equilateral triangles in $\T_n$ contained in a compact subset of $\Omega$ are mapped by $f_n^{-1}$ to triangles in $(\triangle ABC, \T_n)$ which are close to be equilateral. Then $f_n$ restricted in each triangle converge to a similarity map. The dilatations $K_n$ of $f_n$ converge to $1$. Therefore, $g$ is $1$-conformal, which is equivalent to be conformal.
\end{proof}






\end{document}